\documentclass[12pt]{article}

\usepackage{epsfig}
\usepackage{amsmath}
\usepackage{soul,color}
\usepackage{bm}
\usepackage{epsfig}
\usepackage{natbib}
\usepackage{latexsym, graphicx, alltt}
\usepackage{graphicx}
\usepackage{setspace}
\usepackage{hyperref}
\usepackage[table]{xcolor}
\usepackage{float}
\usepackage[nolists]{endfloat}

\def\b#1{\mbox{\boldmath $#1$}}

\newcommand{\cb}{\cellcolor{blue!25}}

\newcommand{\E}{{\rm E}}
\newcommand{\N}{{\rm N}}

\def\cg#1{\mbox{${\cal #1}$}}
\newcommand{\tr}{^{\prime}}

\newcommand{\be}{\beta}

\newcommand{\si}{\sigma}

\newcommand{\la}{\lambda}
\renewcommand{\th}{\theta}

\def\bl#1{\mbox{\scriptsize\boldmath {$#1$}}} 

\newcommand\bv{\bar{v}}
\begin{document}

\title{\vspace*{-1.5cm}Composite likelihood inference in a discrete latent variable model for two-way ``clustering-by-segmentation" problems}

\author{Francesco Bartolucci \\ {\small Department of Economics} \\ 
{\small University of Perugia (IT)} \\ 
{\small email: francesco.bartolucci@unipg.it}
\and
Francesca Chiaromonte \\ {\small Department of Statistics} \\
{\small The Pennsylvania State University (USA)} \\ 
{\small email: 
fxc11@psu.edu} 
\and
Prabhani Kuruppumullage Don \\ 
{\small Department of Biostatistics and Computational Biology} \\
{\small Dana-Farber Cancer Institute, and} \\ 
{\small Department of Biostatistics} \\
{\small Harvard School of Public Health (USA)} \\
{\small email: pdon@jimmy.harvard.edu} \\
\and
Bruce G. Lindsay \\ 
{\small Department of Statistics} \\
{\small The Pennsylvania State University (USA)} \\ 
{\small email: 
bgl@psu.edu}}
\maketitle

\begin{abstract}
We consider a discrete latent variable model for two-way data arrays, which allows one to simultaneously produce clusters along one of the data dimensions (e.g.~exchangeable observational units or features) and contiguous groups, or segments, along the other (e.g.~consecutively ordered times or locations). The model relies on a hidden Markov structure but, given its complexity, cannot be estimated by full maximum likelihood. We therefore introduce composite likelihood methodology based on considering different subsets of the data. The proposed approach is illustrated by simulation, and with an application to genomic data.\vskip3mm

\noindent{\sc Key Words}: Crossed-effects models; Cross validation; EM algorithm; Finite mixture models; Composite likelihood; Genomics.
\end{abstract}\newpage


\section{Introduction}\label{sec:intro}

Many recent applications involve large two-way data arrays in which both rows and columns need to be grouped, possibly taking into account a serial dependence in one of the two dimensions. Applications of this type arise in several fields. For instance, in Economics, they may concern parallel time-series for a certain indicator recorded on a pool of countries. In this case, one may be interested in clustering countries and simultaneously grouping contiguous time periods into segments corresponding to different phases of the economic cycle. As another example, Genomics data sets often comprise a number of features measured along the nuclear DNA of a species, capturing characteristics of the DNA sequence and/or various types of molecular activities. In these settings, one may be interested in clustering such features and simultaneously partitioning the genome into segments corresponding to different molecular activity landscapes. 

To deal with this type of ``clustering-by-segmentation" problems, we introduce a statistical model based on associating a discrete latent variable to every row and column of a two-way data array. The row latent variables are assumed to be independent and identically distributed, as they refer to entities that are exchangeable in nature (e.g.~the countries, or the genomic features). The column latent variables, on the other hand, refer to serially dependent entities (e.g.~time periods, or locations along the nuclear DNA) and are assumed to follow a first-order homogenous hidden Markov (HM) model \citep{zucc:macd:09} with initial distribution equal to the stationary distribution. Given row and column latent variables, the observable variables are assumed to be conditionally independent and distributed according to laws whose parameters depend on the values of the latent variables themselves. Our approach is not restricted to a specific type of outcomes, so that a generalized linear parametrization as in  \cite{mcc:nel:89} may be used to relate observable and latent variables.

Our {\em two-way discrete latent variable model} comprises well-known models as special cases. In particular, it includes the class of crossed random effect models considered by \cite{bellio:varin:2005} when these random effects are assumed to have a discrete distribution. This class of models, however, does not comprise any serial dependence for the column latent variables.

The main focus of this article is likelihood-based inference for our model. As we will show, when the dimensions of the two-way data array are small, maximizing the full model likelihood is computationally feasible and may be performed by an Expectation-Maximization (EM) algorithm \citep{baum:et:al:70,demp:lair:rubi:77} that extends the one used for HM models of longitudinal data \citep{bart:farc:penn:12,bartolucci2014latent} -- also named latent Markov models. 
However, it is easy to convince oneself that full model likelihood maximization is computationally unaffordable for moderate or large size arrays. 
In fact, even employing the efficient and well-known HM forward recursion by Baum and Welch \citep{baum:et:al:70,welch:2003}, the numerical complexity of computing the full likelihood function for an array of dimensions $r\times s$ has order $O(sk_1^rk_2^2)$ -- where $k_1$ and $k_2$ are the number of support points of row and column latent variables, respectively. This complexity increases exponentially with $r$ and linearly with $s$, due the use of the aforementioned HM recursion. Notably then, data arrays with a large number of exchangeable rows are more problematic to deal with than those with a large number of serially dependent columns.

In order to deal with the estimation problem described above, we propose a composite likelihood approach \citep{lindsay1988composite,cox2004note}; see \citep{varin2011overview} for a review.
In particular, we introduce two versions of composite likelihood. The first, which we name {\em row composite likelihood}, results from ignoring dependencies between data rows due to sharing the same column latent variables. This composite likelihood can be maximized by an EM algorithm similar to the one used for mixed HM models of longitudinal data with discrete mixing distributions \citep{maruotti2011mixed}. The second and more satisfactory version, which we name {\em row-column composite likelihood}, results from combining the row composite likelihood with an analogous construct for the columns; i.e.~a composite likelihood in which one ignores dependencies between data columns due to sharing row latent variables.
As we will show, also the row-column composite likelihood can be maximized by an EM algorithm that is computationally viable even for large data arrays. Our algorithms are implemented in {\tt R} and available upon request.

We study the finite sample properties of row and row-column composite likelihood estimators by simulation. Importantly, our simulation study covers also two-way arrays with small dimensions -- where these estimators can in fact be compared with the full likelihood estimator.
This gives us a chance to quantify the loss of efficiency due to the use of composite likelihood approximations, and to identify the parameters with respect to which this loss is more sizable.

Another relevant aspect we tackle is model selection; in particular, the choice of $k_1$ and $k_2$ -- the number of support points for row and column latent variables. The strategy we suggest is based on cross validation -- extending the approach of \cite{smyth2000model} for finite mixture models and of \cite{celeux:durand:2008} for HM models. For our model, implementing cross validation-based model selection is complicated by the lack of independence between any pair of observations in the data. To deal with this, we devise a cross validation scheme where (a) half of the ``cells'' in the two-way data array, identified randomly drawing row and column indexes, are withdrawn for use as test set, and (b) a missing-at-random version of the composite likelihood is used both for estimating parameters on the training set and for measuring fit on the test set.

The use of our model, inference approach and model selection strategy are illustrated through an application in Genomics. Several recent studies \citep{old:et:al:2010,ernst:et:al:2011,encode:2012,hoffman:et:al:2013} have utilized HM models to create segmentations of the human genome leveraging data from inter- or intra-species comparisons, or various types of high-throughput genomic assays. In particular, \cite{kuru:et:al:2013} produced a segmentation based on the rates of four types of mutations estimated from primate comparisons in 1Mb (megabase) non-overlapping windows along the human genome. The authors also gathered and pre-processed publicly available data on several dozens genomic features in the same windows system. These features capture, among other things, aspects of DNA composition, prevalence of transposable elements, recombination rates, chromatin structure, methylation, transcription, etc. Producing a segmentation based on this large array of features could provide significant biological insights, all the more if one could simultaneously characterize their interdependencies by partitioning them into meaningful groups. The application we present in this article is a feasibility proof for such an endeavor; we utilize our model and methodology to perform ``clustering-by-segmentation" on a two-way data array comprising  $r=28$ genomic features measured in $s=224$ contiguous 1Mb windows along human chromosome 1. 
 
The reminder of the article is organized as follows. In Section~\ref{sec:assumptions} we introduce the structure and assumptions underlying our statistical model. In Section~\ref{sec:full_likelihood} we outline methodology for full likelihood estimation of the model parameters, and in Section~\ref{sec:composite_likelihood} we outline row and row-column composite likelihood methodology. In Section~\ref{sec:simulations} we describe our simulation study, in Section~\ref{sec:model selection} we discuss model selection with cross validation, and in Section \ref{sec:illustration} we present results of our application to genomic data. Finally, we offer some concluding remarks in Section~\ref{sec:conclusions}.


\section{The Model}\label{sec:assumptions}

Consider a two-way array of random variables $Y_{ij}$, $i=1,\ldots,r$, $j=1,\ldots,s$, where $r$ is the number of rows and $s$ is the number of columns. The basic assumption of our model is that these observable variables are conditionally independent given two vectors $U_1,\ldots,U_r$ and $V_1,\ldots,V_s$ of {\em row} and  {\em column latent variables}.
The row latent variables ($U$'s) are assumed to be independent and identically distributed according to a discrete distribution with $u=1,\ldots,k_1$ support points and mass probabilities
\[
\la_u=p(U_i=u),\quad 
u=1,\ldots,k_1.
\]
The column latent variables ($V$'s) are assumed to follow a first order Markov chain with $v=1,\ldots,k_2$ states, initial probabilities
\[
\pi_v = p(V_1=v),\quad v=1,\ldots,k_2,
\]
transition probabilities
\[
\pi_{\bv v} = p(V_j=v|V_{j-1}=\bv),\quad 
\bv,v = 1,\ldots,k_2,
\]
and stationary probabilities
\[
\rho_v = \lim_{s\rightarrow\infty} p(V_j=v),\quad v=1,\ldots,k_2.
\]
We also postulate that initial and stationary distributions coincide; that is
\begin{equation}
\pi_v = \rho_v,\quad v=1,\ldots,k_2.\label{eq:stationary}
\end{equation}
This makes the model more parsimonious as the chain can be directly parametrized by the transition probabilities.

Our model specification is completed by formulating the conditional distribution of every observable variable $Y_{ij}$ given the underlying pair of latent variables $(U_i,V_j$).
In the continuous case, a natural assumption is that
\begin{equation}
Y_{ij}|U_i=u,V_j=v\sim \N(\psi_{uv},\si^2),\quad 
u=1,\ldots,k_1,\:v=1,\ldots,k_2,\label{eq:normal}
\end{equation}
where the $\psi_{uv}$ are means depending on the latent variables and $\si^2$ is a common variance. This results in a complex finite mixture of Normal distributions \citep{lindsay1995mixture,mclachlan2004finite}. Note that the requirement that the Normal mixture be homoschedastic
is quite common in the finite mixture literature as it avoids degenerate solutions in terms of maximum likelihood estimates. Moreover in many practical applications (see for instance Section~\ref{sec:illustration}) the data can be preprocessed and transformed as to make a homoschedastic Normal mixture suitable.  

It is important to remark that our model can be made more parsimonious incorporating knowledge in the form of constraints imposed on the means $\psi_{uv}$. For instance, we could postulate that
\[
\psi_{uv}=\psi^{(1)}_u+\psi_v^{(2)},\quad u=1,\ldots,k_1,\:v=1,\ldots,k_2.
\]
On the other hand, our model can be made more general allowing each $Y_{ij}$ to depend also on observable covariates. For instance, we could postulate that
\begin{equation}
\E(Y_{ij}|U_i=u,V_j=v)=\psi_{uv}+\b x_{ij}\tr\b\be,\quad u=1,\ldots,k_1,\:v=1,\ldots,k_2,
\label{eq:norm_cov}
\end{equation}
where the vector $\b x_{ij}$ comprises the covariates (which are assumed to be fixed and known; not random) and the vector $\b\be$ the corresponding regression coefficients (which are assumed not to depend on the latent variables).

Along the same lines adopted in \cite{bellio:varin:2005}, another obvious way to generalize our model is to replace the Normal specification (\ref{eq:normal}) for the distribution of $Y_{ij}$ given $(U_i,V_j)$ with any exponential family distribution used in Generalized Linear Models \citep[GLMs,][]{mcc:nel:89}. For instance, for a two-way array of binary variables we may assume 
\[
Y_{ij}|U_i=u,V_j=v\sim{\rm Bernoulli}(p_{uv}),\quad 
u=1,\ldots,k_1,\:v=1,\ldots,k_2,
\]
where the success probabilities $p_{uv}$ depend on the latent variables. Then, as in a GLM, these probabilities could be expressed through an additive parametrization and/or as a function of observable covariates. For instance, depending on the application, it may be reasonable to postulate that
\[
\log\frac{\E(Y_{ij}=1|U_i=u,V_j=v)}{\E(Y_{ij}=0|U_i=u,V_j=v)}=\psi_{u}^{(1)}+\psi_{v}^{(2)}+\b x_{ij}\tr\b\be,\quad u=1,\ldots,k_1,\:v=1,\ldots,k_v,
\]
using a logit link function -- which directly compares with (\ref{eq:norm_cov}).

In the following, we first introduce full likelihood maximization as a means to estimate the parameters of our model. As this type of estimation is computationally feasible only with small arrays, we then switch to composite likelihood methodology. We remark that while the methods are described in reference to the homoschedastic Normal mixture in (\ref{eq:normal}) and under the constraint that the initial and stationary distributions of the Markov chain coincide as in (\ref{eq:stationary}), the implementation can be streighforwardly generalized to deal with different parameterizations and/or to account for covariates.


\section{Full Likelihood Methodology}\label{sec:full_likelihood}

Let $\b Y$ denote the matrix of all the outcomes $y_{ij}$, $i=1,\ldots,r$, $j=1,\ldots,s$. Also let $\b u=(u_1,\ldots,u_r)\tr$ and $\b v=(v_1,\ldots,v_s)\tr$ denote possible configurations of the row and column latent variables, respectively. The joint density function will then be
\[
p(\b Y)=\sum_{\bl u}\sum_{\bl v}\la_{u_1}\cdots\la_{u_r}\rho_{v_1}\pi_{v_1v_2}
\cdots\pi_{v_{s-1}v_s}\prod_i\prod_j\phi(y_{ij};\psi_{u_iv_j},\si^2),
\]
where $\phi(y;\psi,\si^2)$ denotes the density function of a $\N(\psi,\si^2)$, and the sums 
are extended to all possible row and column latent variables configurations $\b u$ and $\b v$. Expressed this way, $p(\b Y)$ can be computed only in trivial cases because it involves a sum over $k_1^rk_2^s$ terms. However, if the number of rows $r$ is relatively small, an effective strategy is to rewrite the joint density function as
\[
p(\b Y)=\sum_{\bl u}p(\b Y|\b u)p(\b u),
\]
where 
\begin{equation}
p(\b Y|\b u)=\sum_{v_1}\rho_{v_1}p(\b y_1^{(2)}|\b u,v_1)
\sum_{v_2}\pi_{v_1v_2}p(\b y_2^{(2)}|\b u,v_2)
\cdots\sum_{v_s}\pi_{v_{s-1}v_s}p(\b y_s^{(2)}|\b u,v_s),
\label{pcond}
\end{equation}
and $\b y_j^{(2)}=(y_{1j},\ldots,y_{rj})\tr$ corresponds to a single column of outcomes, so that
\[
p(\b y_j^{(2)}|\b u,v)=\prod_j\phi(y_{ij};\psi_{u_iv},\si^2).
\]
The density function in (\ref{pcond}) can be computed with a well-known recursion in the HM literature \citep{baum:et:al:70,welch:2003}, the numerical complexity of which increases linearly in $s$. Thus, the number of operations to compute $p(\b Y)$ becomes of order $sk_1^rk_2$, as already indicated in Section~\ref{sec:intro}. Relatedly, the log-likelihood
\[
\ell(\b\th) = \log p(\b Y),
\]
where $\b\th$ is short-hand notation for all model parameters, can be maximized using an EM algorithm \citep{baum:et:al:70,demp:lair:rubi:77} which is described in detail in the following.


\subsection{EM algorithm for full likelihood estimation}\label{sec:EM_full}

First, we introduce the {\em complete data log-likelihood} corresponding to $\ell(\b\th)$. Consider the latent indicators $w_{iu}$ and $z_{jv}$ -- for row and column latent variables, respectively. In particular, $w_{iu}$ is equal to 1 if $U_i=u$ and to 0 otherwise, with $i=1,\ldots,r$, $u=1,\ldots,k_1$, and $z_{jv}$ is similarly defined with reference to $V_j$. With some simple algebra, the complete data log-likelihood can be written as
\begin{equation}
\ell^*(\b\th) = a(\b\la)+b(\b\Pi)+c(\b\Psi,\si^2),\label{eq:comp_lk}
\end{equation}
where
\begin{eqnarray*}
a(\b\la)&=&\sum_i \sum_u w_{iu} \log(\la_u),\\
b(\b\Pi)&=&\sum_v z_{1v}\log(\rho_v)+\sum_{j>1}\sum_{\bv}\sum_v z_{j\bv v}\log(\pi_{\bv v}),\\
c(\b\Psi,\si^2)&=&\sum_i\sum_j\sum_u\sum_vw_{iu}z_{jv}\log\phi(y_{ij};\psi_{uv},\si^2),\nonumber
\end{eqnarray*}
with $z_{j\bv v}=z_{j\bv}z_{jv}$. In the above decomposition, the vector $\b\la$ comprises the row latent variables' mass probabilities $\la_v$, the matrix $\b\Pi$ comprises the column latent variables' transition probabilities $\pi_{\bv v}$, 
and the matrix $\b\Psi$ comprises the means $\psi_{uv}$. 

The EM algorithm alternates two steps until convergence:
\begin{itemize}
\item {\bf E-step}: compute the posterior expected value of each indicator variable in (\ref{eq:comp_lk}). For $i=1,\ldots,n$ and $u=1,\ldots,k_1$ we set
\[
\hat{w}_{iu} = p(U_i=u|\b Y)=\frac{1}{p(\b Y)}\sum_{\bl u:u_i=u}p(\b Y|\b u)p(\b u),
\]
where the sum $\sum_{\bl u:u_i=u}$ is extended to all configurations $\b u$ with $i$th element equal to $u$. For $j=1,\ldots,s$ and $\bv,v=1,\ldots,k_2$ we set
\begin{eqnarray*}
\hat{z}_{1v}&=& p(V_1=v|\b Y)=\frac{1}{p(\b Y)}\sum_{\bl u}p(V_1=v|\b u,\b Y)p(\b Y|\b u)p(\b u),\\
\hat{z}_{j\bv v}&=& p(V_{j-1}=\bv,V_j=v|\b Y)=
\frac{1}{p(\b Y)}\sum_{\bl u}p(V_{j-1}=\bv,V_j=v|\b u,\b Y)p(\b Y|\b u)p(\b u),
\end{eqnarray*}
where the conditional probabilities $p(V_j=v|\b u,\b Y)$ and $p(V_{j-1}=\bv,V_j=v|\b u,\b Y)$ are obtained from suitable recursions \citep{baum:et:al:70,welch:2003}. Finally, for $i=1,\ldots,n$, $j=1,\ldots,s$, $u=1,\ldots,k_1$, and $\bv,v=1,\ldots,k_2$ we set
\[
\widehat{(w_{iu}z_{jv})}=p(U_i=u,V_j=v|\b Y)=
\frac{1}{p(\b Y)}\sum_{\bl u:u_i=u}p(V_j=v|\b u,\b Y)p(\b Y|\b u)p(\b u).
\]
\item {\bf M-step}: update the value of each parameter in (\ref{eq:comp_lk}). For $u=1,\ldots,k_1$ we update the row mass probabilities as
\[
\la_u = \frac{1}{r}\sum_i \hat{w}_{iu}.
\]
Under constraint (\ref{eq:stationary}), we update the transition probabilities by numerical maximization of the function
\[
\hat{b}(\b\Pi)=\sum_v \hat{z}_{1v}\log(\rho_v)+\sum_{j>1}\sum_{\bv}\sum_v \widehat{z}_{j\bv v}\log(\pi_{\bv v});
\]
see also \cite{bul-ber:2008} and \cite{zucc:macd:09}.
Finally, for $u=1,\ldots,k_1$ and $v=1,\ldots,k_2$ we update means and common variance for the Normal distributions as 
\[
\mu_{uv}=\frac{1}{\sum_i\sum_j\widehat{(w_{iu}z_{jv})}}
\sum_i\sum_j\widehat{(w_{iu}z_{jv})}y_{ij}
\]
and
\[
\si^2=\frac{1}{rs}
\sum_i\sum_j\sum_u\sum_v\widehat{(w_{iu}z_{jv})}(y_{ij}-\mu_{uv})^2.
\]
\end{itemize}
This algorithm runs and converges in a reasonable time if the number of rows in the two-way data array is $r \leq 10$ and the row latent variables are binary ($k_1=2$), even with a large number $s$ of columns. Just to give an idea, using our {\tt R} implementation on a standard personal computer, a few seconds are necessary to estimate the model with $r=5$ rows, $s=200$ columns, and binary latent variables. Again with binary latent variables and the same value of $s$, but with $r=10$, the computing time increases to a few minutes.
However, as $r$ increases and, in particular, as the number of support points of the row latent variables increases, full maximum likelihood estimation becomes prohibitive and is infeasible for the models considered in our application (see Section \ref{sec:illustration}). 
\section{Composite Likelihood Methodology}\label{sec:composite_likelihood}
Given that the EM algorithm for full likelihood estimation is not computationally viable in typical applications, we propose an alternative approach based on maximizing a composite likelihood function where the rows are treated separately \citep{lindsay1988composite,cox2004note}. In this section we introduce two versions of the composite likelihood function; the row composite likelihood, which is related to the method proposed by \citep{bartolucci2015pairwise} for multilevel HM models, and the row-column composite likelihood. The latter is characterized by greater complexity and potentially larger estimation efficiency.
\subsection{Row composite likelihood estimation}

First, we consider the density function of the $i$th row of the data, represented as a column vector $\b y_i^{(1)}=(y_{i1},\ldots,y_{is})\tr$. Given the underlying latent variable $U_i$, this is generated along a stationary hidden Markov model, so that
\[
p(\b y^{(1)}_i|U_i=u)=\sum_{v_1}\rho_{v_1}\phi(y_{i1};\psi_{uv_1},\si^2)
\sum_{v_2} \pi_{v_1v_2} \phi(y_{i2};\psi_{uv_2},\si^2)\cdots\sum_{v_s} 
\pi_{v_{s-1}v_s}\phi(y_{is};\psi_{uv_s},\si^2).
\] 
In practice, $p(\b y^{(1)}_i|U_i=u)$ is computed by a simplified version of the recursion used for the full likelihood estimation. The next step is to integrate out the latent variable $U_i$ as to obtain
\[
p(\b y^{(1)}_i)=\sum_u\la_up(\b y^{(1)}_i|U_i=u).
\]
The {\em row composite log-likelihood} is defined based on this density function as
\begin{equation}
c\ell_1(\b\th) = \sum_i \log p(\b y^{(1)}_i).\label{eq:composite_lk}
\end{equation}
Importantly, this can be readily computed also for a large number of rows, as it treats the rows as independent.

In order to implement an EM algorithm to maximize $c\ell(\b\th)$, it is useful to note that (\ref{eq:composite_lk}) is the log-likelihood of a model that, in addition to satisfying all assumptions in Section~\ref{sec:assumptions}, postulates independent Markov chains $V_{i1},\ldots,V_{is}$ underlying each row of data $\b y_i^{(1)}$, $i=1,\ldots,r$. This additional assumption implies a different definition of the complete data likelihood. We now need to consider the indicator variables $w^{(1)}_{iu}$ and $z^{(1)}_{ijv}$. 
The former have the same meaning as the $w_{iu}$ introduced in Section~\ref{sec:assumptions}, however the latter are now defined separately for each row -- reflecting the structure of the target function in (\ref{eq:comp_lk}); we let $z^{(1)}_{ijv}$ equal to 1 if $V_{ij}=v$ and to 0 otherwise. Using these indicator variables, we express the complete data composite log-likelihood as
\begin{equation}
c\ell_1^*(\b\th) = ca_1(\b\la)+cb_1(\b\Pi)+cc_1(\b\Psi,\si^2),\label{eq:composite_comp_lk}
\end{equation}
where
\begin{eqnarray*}
ca_1(\b\la)&=&\sum_i \sum_u w_{iu}^{(1)} \log(\la_u),\\
cb_1(\b\Pi)&=&\sum_i\left[\sum_v z_{i1v}^{(1)}\log(\rho_v)+
\sum_{j>1}\sum_{\bv}\sum_v z_{ij\bv v}^{(1)}\log(\pi_{\bv v})\right],\\
cc_1(\b\Psi,\si^2)&=&\sum_i\sum_j\sum_u\sum_vw_{iu}^{(1)}z_{ijv}^{(1)}\log\phi(y_{ij};\psi_{uv},\si^2),\nonumber
\end{eqnarray*}
with $z_{ij\bv v}^{(1)}=z_{i,j-1,\bv v}^{(1)}z_{ijv}^{(1)}$.

The EM alternates two steps until convergence:
\begin{itemize}%
\item {\bf E-step}: compute the  posterior expected value of each indicator variable in (\ref{eq:composite_comp_lk}). Note the definitions in terms of posterior probabilities here hold under the ``approximating'' model in which the data rows are independent. For $i=1,\ldots,r$ and $u=1,\ldots,k_1$ we set
\begin{equation}
\hat{w}_{iu}^{(1)} = p(U_i=u|\b y_i^{(1)})=\frac{p(\b y_i^{(1)}|U_i=u)\la_u}{p(\b y_i^{(1)})}.
\label{eq:ppost11}
\end{equation}
Thus, for $i=1,\ldots,r$, $j=2,\ldots,s$, and $\bv,v=1,\ldots,k_2$ we set
\begin{eqnarray}
\hat{z}^{(1)}_{i1v}&=&p(V_{i1}=v|\b y_i^{(1)})=
\sum_up(V_{i1}=v|U_i=u,\b y_i^{(1)})\hat{w}_{iu}^{(1)},\label{eq:ppost12}\\
\hat{z}^{(1)}_{ij\bv v}&=& p(V_{i,j-1}=\bv,V_{ij}=v|\b y_i^{(1)})=
\sum_up(V_{i,j-1}=\bv,V_{ij}=v|U_i=u,\b y_i^{(1)})\hat{w}_{iu}^{(1)},\label{eq:ppost13}
\end{eqnarray}
where the conditional probabilities $p(V_{ij}=v|U_i=u,\b y_i^{(1)})$ and $p(V_{ij}=v,V_{i,j-1}=\bv|U_i=u,\b y_i^{(1)})$ are obtained by suitable recursions similar to the ones used in the E-step for the full likelihood. Finally, for $=1,\ldots,r$, $j=1,\ldots,s$, $u=1,\ldots,k_1$ and $v=1,\ldots,k_2$ we set
\begin{equation}
\widehat{(w^{(1)}_{iu}z^{(1)}_{ijv})}=p(U_i=u,V_{ij}=v|\b y_i^{(1)})=
p(V_{ij}=v|U_i=u,\b y_i^{(1)})\hat{w}_{iu}^{(1)}.\label{eq:ppost14}
\end{equation}
\item {\bf M-step}: update the value of each parameter in (\ref{eq:composite_comp_lk}). For $u=1,\ldots,k_1$ we update the row mass probabilities as
\[
\la_u = \frac{1}{r}\sum_i \hat{w}_{iu}^{(1)}. 
\]
Under constraint (\ref{eq:stationary}), we update the transition probabilities by numerical maximization of the function
\[
\widehat{cb}_1(\b\Pi)=\sum_i\left[\sum_u\sum_v \hat{z}^{(1)}_{ijv}\log(\rho_v)+\sum_{j>1}\sum_{\bv}\sum_v \hat{z}^{(1)}_{ij\bv v}\log(\pi_{\bv v})\right].
\]
Finally, for $u=1,\ldots,k_1$ and $v=1,\ldots,k_2$ we update means and common variance for the Normal distributions as
\[
\mu_{uv}=\frac{\sum_i\sum_j\widehat{(w^{(1)}_{iu}z^{(1)}_{ijv})}y_{ij}}
{\sum_i\sum_j\widehat{(w^{(1)}_{iu}z^{(1)}_{ijv})}}
\]
and
\[
\si^2=\frac{1}{rs}
\sum_i\sum_j\sum_u\sum_v\widehat{(w^{(1)}_{iu}z^{(1)}_{ijv})}(y_{ij}-\mu_{uv})^2.
\]
\end{itemize}
%


\subsection{Row-column composite likelihood estimation}

We now pass to consider a more complex composite likelihood, which takes into account also the density function of each separate column of the data. For the $j$th data column represented by $\b y_j^{(2)}$, given the underlying latent variable $V_j$, we have
\[
p(\b y_j^{(2)}|V_j=v)=\prod_ip(y_{ij}|V_j=v),
\] 
where $p(y_{ij}|V_j=v)=\sum_u \phi(y_{ij};\psi_{uv},\si^2)\la_u$. Thus, integrating out the latent variable $V_j$ we obtain
\[
p(\b y^{(2)}_j)=\sum_vp(\b y^{(2)}_j|V_j=v)\rho_v.
\]
The composite log-likelihood based on this density function is
\begin{equation}
c\ell_2(\b\th) = \sum_j \log p(\b y^{(2)}_j).\label{eq:composite_lk2}
\end{equation}
To estimate the parameters of our model, we propose to maximize the {\em row-column composite log-likelihood} defined as the sum of the row composite log-likelihood in (\ref{eq:composite_lk}) with the above expression:
\[
c\ell(\b\th) = c\ell_1(\b\th)+c\ell_2(\b\th)
\]
In this regard, we note that (\ref{eq:composite_lk2}) is the log-likelihood of a model which, in addition to satisfying all assumptions in Section~\ref{sec:assumptions}, postulates that each column of the data $\b y_j^{(2)}$, $j=1,\ldots,s$, depends on an independent sequence of latent variables $U_{1j}, \ldots U_{rj}$ also independent of each other. Moreover, this model assumes that the latent variables $V_j$, $j=1,\ldots,s$ are independent and distributed according to the stationary distribution. Consequently, we now use the indicator variables $w^{(2)}_{iju}$ and $z^{(2)}_{jv}$. 
The latter have the same meaning as the $z_{jv}$ introduced in Section~\ref{sec:EM_full}, however the former are defined separately for each column; we set $w^{(2)}_{iju}=1$ if $U_{ij}=u$ and $0$ otherwise. Using these indicator variables, we express the complete data composite log-likelihood as
\[
c\ell_2^*(\b\th) = ca_2(\b\la)+cb_2(\b\Pi)+cc_2(\b\Psi,\si^2),
\]
where
\begin{eqnarray*}
ca_2(\b\la)&=&\sum_i\sum_j\sum_u w_{iju}^{(2)} \log(\la_u),\\
cb_2(\b\Pi)&=&\sum_j\sum_v z_{jv}^{(2)}\log(\rho_v),\\
cc_2(\b\Psi,\si^2)&=&\sum_i\sum_j\sum_u\sum_vw_{iju}^{(2)}z_{jv}^{(2)}\log\phi(y_{ij};\psi_{uv},\si^2).\nonumber
\end{eqnarray*}

The EM alternates two steps until convergence:
\begin{itemize}%
\item {\bf E-step}: We compute the same posterior probabilities as in (\ref{eq:ppost11}), (\ref{eq:ppost12}), and (\ref{eq:ppost13}). In addition, for $j=1,\ldots,s$ and $v=1,\ldots,k_2$ we set  
\[
\hat{z}_{jv}^{(2)} = p(V_j=v|\b y_j^{(2)})=\frac{p(\b y_j^{(2)}|V_j=v)\rho_v}
{p(\b y_j^{(2)})}.
\]
Thus, for $i=1,\ldots,r$ and $j=1,\ldots,s,\:u=1,\ldots,k_1$ we set
\[
\hat{w}_{iju}^{(2)} = p(U_{ij}=u|\b y_j^{(2)})=\sum_v \frac{\phi(y_{ij};\psi_{uv},\si^2)\la_u}{p(y_{ij}|V_j=v)}\hat{z}_{jv}^{(2)} 
\]
and for $i=1,\ldots,r$, $j=1,\ldots,s$, $u=1,\ldots,k_1$ and  $v=1,\ldots,k_2$ we set
\[
\widehat{(w_{iu}^{(2)}\hat{z}_{jv}^{(2)})}=p(U_{ij}=u,V_j=v|\b y_j^{(2)})\frac{\phi(y_{ij};\psi_{uv},\si^2)\la_u}{p(y_{ij}|V_j=v)}\hat{z}_{jv}^{(2)}.
\]
\item {\bf M-step}: For $u=1,\ldots,k_1$ we update the row mass probabilities as
\begin{equation}
\la_u = \frac{1}{r+rs}\left[\sum_i \hat{w}_{iu}^{(1)}+\sum_i\sum_j \hat{w}_{iju}^{(2)}\right].\label{eq:update_la_composite}
\end{equation}
We update the transition probabilities by numerical maximization of the function
\[
\widehat{cb}(\b\Pi)=\widehat{cb}_1(\b\Pi)+\widehat{cb}_2(\b\Pi),
\]
where
\[
\widehat{cb}_2(\b\Pi)=\sum_j\sum_v \hat{z}^{(2)}_{jv}\log(\rho_v).
\]
Finally, for $u=1,\ldots,k_1$ and $v=1,\ldots,k_2$ we update the means and common variance of the Normal distributions as
\[
\mu_{uv}=\frac{\sum_i\sum_j[\widehat{(w^{(1)}_{iu}z^{(1)}_{ijv})}+\widehat{(w^{(2)}_{iu}z^{(2)}_{ijv})}]y_{ij}}
{\sum_i\sum_j\widehat{(w^{(1)}_{iu}z^{(1)}_{ijv})}+\widehat{(w^{(2)}_{iu}z^{(2)}_{ijv})}}
\]
and
\[
\si^2=\frac{1}{rs}
\sum_i\sum_j\sum_u\sum_v[\widehat{(w^{(1)}_{iu}z^{(1)}_{ijv})}+\widehat{(w^{(2)}_{iu}z^{(2)}_{ijv})}](y_{ij}-\mu_{uv})^2. 
\]
\end{itemize}


\section{Simulation study}\label{sec:simulations}

We performed 
a simulation study to assess and compare the performance of our two approximations -- the row and the row-column composite likelihoods -- to one another and to full likelihood estimation.

\subsection{Simulation design}
We consider a {\em benchmark design} in which the two-way data array has dimensions 
$r=10$ by $s=200$, with two support points for both row and column latent variables ($k_1=k_2=2$).
This design has $r<<s$, as is perhaps typical in many applications, and is small enough for full likelihood estimation to be viable. We fix the model parameters as follows:
\begin{itemize}
\item $\b\la = (0.5,0.5)\tr$;
\item $\b\Pi = \begin{pmatrix}
0.8808 & 0.1192\\
0.1192 & 0.8808 
\end{pmatrix}$, so that 
$\b\rho = (0.5,0.5)\tr$;
\item $\b\Psi = \begin{pmatrix}
1 & 2\\
3 & 4
\end{pmatrix}$;
\item $\si^2=0.5$.
\end{itemize}
In order to assess the behavior of the estimators under comparison, we also consider other scenarios in which specific elements of the benchmark design are suitably modified. In particular, we consider the following scenarios:
\begin{itemize}
\item $r=15$ instead of $r=10$ and parameters fixed as above; 
\item $s=400$ instead of $s=200$ and parameters fixed as above; 
\end{itemize}
in these scenarios there is a larger amount of information on the structure underlying, respectively, the serially dependent columns or the exchangeable rows.
\begin{itemize}
\item $k_1=3$ instead of $k_1=2$ and parameters fixed as above apart from $\b\la = (1/3,1/3,1/3)\tr$; 
\item $k_2=3$ instead of $k_2=2$ and parameters fixed as above apart from $\b\Pi = \begin{pmatrix}
0.7870 & 0.1065 & 0.1065\\
0.1065 & 0.7870 & 0.1065\\
0.1065 & 0.1065 & 0.7870\\
\end{pmatrix}$, and thus 
$\b\rho = (1/3,1/3,1/3)\tr$; 
\end{itemize}
in these scenarios there is a larger complexity of, respectively, the row or column latent structure.
\begin{itemize}
\item $\si^2=1$ instead of $\si^2=0.5$ and parameters fixed as above; 
\end{itemize}
in this scenario there is a smaller separation between latent states. 

\subsection{Simulation results}

Each scenario is simulated 1,000 times independently, and bias and square root of the mean squared error (RMSE) for parameter estimation are computed for each estimation method -- i.e.~full likelihood, row composite likelihood, and row-column composite likelihood. Results for $\la_u$ are reported in Table~\ref{tab:simula_la}, those for $\pi_{\bv v}$ in Table~\ref{tab:simula_Pi}, those for $\psi_{uv}$ in Table~\ref{tab:simula_Psi}, and those for $\si^2$ in Table~\ref{tab:simula_si2}. In Table~\ref{tab:simula_time} we also report median computing times in seconds, along with median absolute deviations -- which are an important elements in comparing estimation methods.

\begin{table}\centering
{\scriptsize
\begin{tabular}{rrrrrlrrrrrrrrrrrrrrrr}\hline\hline
	&		&		&		&		&		&	\multicolumn3c{Full likelihood}					&&	\multicolumn3c{Row comp. lik.}					&&	\multicolumn3c{Row-column comp. lik.}					\\\cline{7-9}\cline{11-13}\cline{15-17}
\multicolumn1c{$r$}	&	\multicolumn1c{$s$}	&	\multicolumn1c{$k_1$}	&	\multicolumn1c{$k_2$}	&	\multicolumn1c{$\si^2$}	&		&	\multicolumn1c{$\la_1$}	&	\multicolumn1c{$\la_2$}	&	\multicolumn1c{$\la_3$}	&&	\multicolumn1c{$\la_1$}	&	\multicolumn1c{$\la_2$}	&	\multicolumn1c{$\la_3$}	&&	\multicolumn1c{$\la_1$}	&	\multicolumn1c{$\la_2$}	&	\multicolumn1c{$\la_3$}	\\\hline
10	&	200	&	2	&	2	&	0.5	&	bias	&	-0.013	&	0.013	&		&&	-0.013	&	0.013	&		&&	-0.013	&	0.013	&		\\
	&		&		&		&		&	rmse	&	0.157	&	0.157	&		&&	0.157	&	0.157	&		&&	0.155	&	0.155	&		\\\hline
15	&	200	&	2	&	2	&	0.5	&	bias	&	-0.002	&	0.002	&		&&	-0.002	&	0.002	&		&&	-0.002	&	0.002	&		\\
	&		&		&		&		&	rmse	&	0.126	&	0.126	&		&&	0.126	&	0.126	&		&&	0.126	&	0.126	&		\\\hline
10	&	400	&	2	&	2	&	0.5	&	bias	&	0.005	&	-0.005	&		&&	0.002	&	-0.002	&		&&	0.003	&	-0.003	&		\\
	&		&		&		&		&	rmse	&	0.149	&	0.149	&		&&	0.145	&	0.145	&		&&	0.145	&	0.145	&		\\\hline
10	&	200	&	3	&	2	&	0.5	&	bias	&	-0.001	&	0.006	&	-0.005	&&	0.001	&	0.003	&	-0.004	&&	0.002	&	0.002	&	-0.003	\\
	&		&		&		&		&	rmse	&	0.137	&	0.140	&	0.136	&&	0.135	&	0.140	&	0.137	&&	0.134	&	0.137	&	0.135	\\\hline
10	&	200	&	2	&	3	&	0.5	&	bias	&	0.009	&	-0.009	&		&&	0.009	&	-0.009	&		&&	0.009	&	-0.009	&		\\
	&		&		&		&		&	rmse	&	0.150	&	0.150	&		&&	0.150	&	0.150	&		&&	0.149	&	0.149	&		\\\hline
10	&	200	&	2	&	2	&	1.0	&	bias	&	-0.007	&	0.007	&		&&	-0.007	&	0.007	&		&&	-0.007	&	0.007	&		\\
	&		&		&		&		&	rmse	&	0.162	&	0.162	&		&&	0.162	&	0.162	&		&&	0.161	&	0.161	&		\\\hline
\end{tabular}}
\caption{Estimation of the $\la_u$ parameters}\label{tab:simula_la}
\end{table}

\begin{table}\centering
{\scriptsize
\begin{tabular}{rrrrrlrrrrrrrrrrrrrrrrr}\hline\hline
	&		&		&		&		&		&		&	\multicolumn3c{Full likelihood}					&&	\multicolumn3c{Row comp. lik.}					&&	\multicolumn3c{Row-column comp. lik.}					\\\cline{8-10}\cline{12-14}\cline{16-18}
\multicolumn1c{$r$}	&	\multicolumn1c{$s$}	&	\multicolumn1c{$k_1$}	&	\multicolumn1c{$k_2$}	&	\multicolumn1c{$\si^2$}	&		&	\multicolumn1c{$\bv$}	&	\multicolumn1c{$\pi_{\bv1}$}	&	\multicolumn1c{$\pi_{\bv2}$}	&	\multicolumn1c{$\pi_{\bv3}$}	&&	\multicolumn1c{$\pi_{\bv1}$}	&	\multicolumn1c{$\pi_{\bv2}$}	&	\multicolumn1c{$\pi_{\bv3}$}	&&	\multicolumn1c{$\pi_{\bv1}$}	&	\multicolumn1c{$\pi_{\bv2}$}	&	\multicolumn1c{$\pi_{\bv3}$}	\\\hline
10	&	200	&	2	&	2	&	0.5	&	bias	&	1	&	-0.004	&	0.004	&		&&	-0.009	&	0.009	&		&&	-0.008	&	0.008	&		\\
	&		&		&		&		&		&	2	&	0.004	&	-0.004	&		&&	0.009	&	-0.009	&		&&	0.008	&	-0.008	&		\\\cline{6-18}
	&		&		&		&		&	rmse	&	1	&	0.034	&	0.034	&		&&	0.046	&	0.046	&		&&	0.042	&	0.042	&		\\
	&		&		&		&		&		&	2	&	0.035	&	0.035	&		&&	0.045	&	0.045	&		&&	0.042	&	0.042	&		\\\hline
15	&	200	&	2	&	2	&	0.5	&	bias	&	1	&	-0.005	&	0.005	&		&&	-0.010	&	0.010	&		&&	-0.009	&	0.009	&		\\
	&		&		&		&		&		&	2	&	0.005	&	-0.005	&		&&	0.009	&	-0.009	&		&&	0.009	&	-0.009	&		\\\cline{6-18}
	&		&		&		&		&	rmse	&	1	&	0.034	&	0.034	&		&&	0.045	&	0.045	&		&&	0.042	&	0.042	&		\\
	&		&		&		&		&		&	2	&	0.034	&	0.034	&		&&	0.044	&	0.044	&		&&	0.041	&	0.041	&		\\\hline
10	&	400	&	2	&	2	&	0.5	&	bias	&	1	&	-0.007	&	0.007	&		&&	-0.009	&	0.009	&		&&	-0.008	&	0.008	&		\\
	&		&		&		&		&		&	2	&	0.001	&	-0.001	&		&&	0.003	&	-0.003	&		&&	0.002	&	-0.002	&		\\\cline{6-18}
	&		&		&		&		&	rmse	&	1	&	0.027	&	0.027	&		&&	0.035	&	0.035	&		&&	0.033	&	0.033	&		\\
	&		&		&		&		&		&	2	&	0.026	&	0.026	&		&&	0.035	&	0.035	&		&&	0.032	&	0.032	&		\\\hline
10	&	200	&	3	&	2	&	0.5	&	bias	&	1	&	-0.007	&	0.007	&		&&	-0.013	&	0.013	&		&&	-0.010	&	0.010	&		\\
	&		&		&		&		&		&	2	&	0.005	&	-0.005	&		&&	0.010	&	-0.010	&		&&	0.008	&	-0.008	&		\\\cline{6-18}
	&		&		&		&		&	rmse	&	1	&	0.038	&	0.038	&		&&	0.052	&	0.052	&		&&	0.048	&	0.048	&		\\
	&		&		&		&		&		&	2	&	0.032	&	0.032	&		&&	0.046	&	0.046	&		&&	0.042	&	0.042	&		\\\hline
10	&	200	&	2	&	3	&	0.5	&	bias	&	1	&	-0.010	&	0.006	&	0.004	&&	-0.008	&	0.004	&	0.004	&&	-0.011	&	0.005	&	0.006	\\
	&		&		&		&		&		&	2	&	0.003	&	-0.009	&	0.005	&&	0.001	&	-0.001	&	0.000	&&	0.005	&	-0.008	&	0.004	\\
	&		&		&		&		&		&	3	&	0.004	&	0.004	&	-0.007	&&	0.002	&	0.004	&	-0.006	&&	0.005	&	0.003	&	-0.008	\\\cline{6-18}
	&		&		&		&		&	rmse	&	1	&	0.057	&	0.044	&	0.042	&&	0.071	&	0.073	&	0.056	&&	0.065	&	0.061	&	0.052	\\
	&		&		&		&		&		&	2	&	0.042	&	0.056	&	0.043	&&	0.069	&	0.082	&	0.065	&&	0.062	&	0.077	&	0.058	\\
	&		&		&		&		&		&	3	&	0.041	&	0.041	&	0.052	&&	0.056	&	0.071	&	0.069	&&	0.052	&	0.060	&	0.061	\\\hline
10	&	200	&	2	&	2	&	1.0	&	bias	&	1	&	-0.004	&	0.004	&		&&	-0.015	&	0.015	&		&&	-0.004	&	0.004	&		\\
	&		&		&		&		&		&	2	&	0.007	&	-0.007	&		&&	0.019	&	-0.019	&		&&	0.007	&	-0.007	&		\\\cline{6-18}
	&		&		&		&		&	rmse	&	1	&	0.038	&	0.038	&		&&	0.061	&	0.061	&		&&	0.048	&	0.048	&		\\
	&		&		&		&		&		&	2	&	0.036	&	0.036	&		&&	0.062	&	0.062	&		&&	0.049	&	0.049	&		\\\hline
\end{tabular}}
\caption{Estimation of the $\pi_{\bv v}$ parameters}\label{tab:simula_Pi}
\end{table}

\begin{table}\centering
{\scriptsize
\begin{tabular}{rrrrrlrrrrrrrrrrrrrrrrr}\hline\hline
	&		&		&		&		&		&		&	\multicolumn3c{Full likelihood}					&&	\multicolumn3c{Row comp. lik.}					&&	\multicolumn3c{Row-column comp. Lik.}					\\\cline{8-10}\cline{12-14}\cline{16-18}
\multicolumn1c{$r$}	&	\multicolumn1c{$s$}	&	\multicolumn1c{$k_1$}	&	\multicolumn1c{$k_2$}	&	\multicolumn1c{$\si^2$}	&		&	\multicolumn1c{$u$}	&	\multicolumn1c{$\psi_{u1}$}	&	\multicolumn1c{$\psi_{u2}$}	&	\multicolumn1c{$\psi_{u3}$}	&&	\multicolumn1c{$\psi_{u1}$}	&	\multicolumn1c{$\psi_{u2}$}	&	\multicolumn1c{$\psi_{u3}$}	&&	\multicolumn1c{$\psi_{u1}$}	&	\multicolumn1c{$\psi_{u2}$}	&	\multicolumn1c{$\psi_{u3}$}	\\\hline
10	&	200	&	2	&	2	&	0.5	&	bias	&	1	&	0.003	&	0.001	&		&&	-0.002	&	0.006	&		&&	-0.005	&	-0.012	&		\\
	&		&		&		&		&		&	2	&	-0.002	&	-0.002	&		&&	-0.005	&	0.001	&		&&	0.010	&	0.002	&		\\\cline{6-18}
	&		&		&		&		&	rmse	&	1	&	0.071	&	0.072	&		&&	0.083	&	0.084	&		&&	0.076	&	0.076	&		\\
	&		&		&		&		&		&	2	&	0.073	&	0.071	&		&&	0.082	&	0.083	&		&&	0.075	&	0.079	&		\\\hline
15	&	200	&	2	&	2	&	0.5	&	bias	&	1	&	0.002	&	0.001	&		&&	-0.003	&	0.004	&		&&	-0.002	&	-0.003	&		\\
	&		&		&		&		&		&	2	&	0.000	&	-0.001	&		&&	-0.004	&	0.002	&		&&	0.004	&	0.002	&		\\\cline{6-18}
	&		&		&		&		&	rmse	&	1	&	0.027	&	0.028	&		&&	0.045	&	0.045	&		&&	0.035	&	0.033	&		\\
	&		&		&		&		&		&	2	&	0.027	&	0.029	&		&&	0.043	&	0.044	&		&&	0.031	&	0.035	&		\\\hline
10	&	400	&	2	&	2	&	0.5	&	bias	&	1	&	-0.003	&	-0.003	&		&&	-0.004	&	-0.002	&		&&	-0.009	&	-0.017	&		\\
	&		&		&		&		&		&	2	&	-0.013	&	-0.014	&		&&	-0.017	&	-0.013	&		&&	0.000	&	-0.009	&		\\\cline{6-18}
	&		&		&		&		&	rmse	&	1	&	0.025	&	0.024	&		&&	0.034	&	0.040	&		&&	0.030	&	0.034	&		\\
	&		&		&		&		&		&	2	&	0.160	&	0.153	&		&&	0.170	&	0.168	&		&&	0.169	&	0.163	&		\\\hline
10	&	200	&	3	&	2	&	0.5	&	bias	&	1	&	0.009	&	0.013	&		&&	0.005	&	0.019	&		&&	-0.001	&	-0.005	&		\\
	&		&		&		&		&		&	2	&	-0.024	&	-0.027	&		&&	-0.029	&	-0.021	&		&&	-0.013	&	-0.036	&		\\
	&		&		&		&		&		&	3	&	-0.051	&	-0.052	&		&&	-0.055	&	-0.047	&		&&	-0.032	&	-0.041	&		\\\cline{6-18}
	&		&		&		&		&	rmse	&	1	&	0.150	&	0.151	&		&&	0.158	&	0.159	&		&&	0.156	&	0.153	&		\\
	&		&		&		&		&		&	2	&	0.304	&	0.302	&		&&	0.302	&	0.303	&		&&	0.303	&	0.304	&		\\
	&		&		&		&		&		&	3	&	0.319	&	0.320	&		&&	0.320	&	0.321	&		&&	0.316	&	0.322	&		\\\hline
10	&	200	&	2	&	3	&	0.5	&	bias	&	1	&	0.007	&	0.007	&	0.005	&&	0.008	&	0.004	&	0.002	&&	0.005	&	-0.002	&	-0.002	\\
	&		&		&		&		&		&	2	&	0.001	&	0.002	&	-0.001	&&	0.004	&	-0.002	&	-0.002	&&	0.007	&	0.012	&	0.001	\\\cline{6-18}
	&		&		&		&		&	rmse	&	1	&	0.144	&	0.143	&	0.133	&&	0.157	&	0.188	&	0.152	&&	0.145	&	0.143	&	0.141	\\
	&		&		&		&		&		&	2	&	0.047	&	0.044	&	0.043	&&	0.081	&	0.141	&	0.083	&&	0.051	&	0.055	&	0.051	\\\hline
10	&	200	&	2	&	2	&	1.0	&	bias	&	1	&	0.003	&	0.003	&		&&	-0.009	&	0.019	&		&&	0.011	&	-0.033	&		\\
	&		&		&		&		&		&	2	&	0.000	&	0.001	&		&&	-0.012	&	0.017	&		&&	0.038	&	-0.004	&		\\\cline{6-18}
	&		&		&		&		&	rmse	&	1	&	0.102	&	0.103	&		&&	0.130	&	0.132	&		&&	0.116	&	0.118	&		\\
	&		&		&		&		&		&	2	&	0.083	&	0.083	&		&&	0.114	&	0.119	&		&&	0.103	&	0.098	&		\\\hline
\end{tabular}}
\caption{Estimation of the $\psi_{uv}$ parameters}\label{tab:simula_Psi}
\end{table}

\begin{table}\centering
{\scriptsize
\begin{tabular}{rrrrrlrrr}\hline\hline
\multicolumn1c{$r$}	&	\multicolumn1c{$s$}	&	\multicolumn1c{$k_1$}	&	\multicolumn1c{$k_2$}	&	\multicolumn1c{$\si^2$}	&		&	\multicolumn1c{Full likelihood}	&	\multicolumn1c{Row likelihood}	&	\multicolumn1c{Row-column likelihood}	\\\hline
10	&	200	&	2	&	2	&	0.5	&	bias	&	-0.001	&	-0.005	&	-0.003	\\
	&		&		&		&		&	rmse	&	0.016	&	0.028	&	0.020	\\\hline
15	&	200	&	2	&	2	&	0.5	&	bias	&	-0.001	&	-0.004	&	-0.003	\\
	&		&		&		&		&	rmse	&	0.013	&	0.025	&	0.016	\\\hline
10	&	400	&	2	&	2	&	0.5	&	bias	&	0.000	&	-0.002	&	-0.002	\\
	&		&		&		&		&	rmse	&	0.012	&	0.021	&	0.015	\\\hline
10	&	200	&	3	&	2	&	0.5	&	bias	&	-0.001	&	-0.006	&	-0.003	\\
	&		&		&		&		&	rmse	&	0.016	&	0.030	&	0.022	\\\hline
10	&	200	&	2	&	3	&	0.5	&	bias	&	-0.001	&	0.002	&	-0.001	\\
	&		&		&		&		&	rmse	&	0.017	&	0.039	&	0.019	\\\hline
10	&	200	&	2	&	2	&	1.0	&	bias	&	-0.002	&	-0.015	&	0.006	\\
	&		&		&		&		&	rmse	&	0.032	&	0.057	&	0.039	\\\hline\end{tabular}}
\caption{Estimation of $\si^2$}\label{tab:simula_si2}
\end{table}

\begin{table}\centering
{\scriptsize
\begin{tabular}{rrrrrrrrrrr}\hline\hline
\multicolumn1c{$r$}	&	\multicolumn1c{$s$}	&	\multicolumn1c{$k_1$}	&	\multicolumn1c{$k_2$}	&	\multicolumn1c{$\si^2$}	&	\multicolumn2c{Full likelihood}	&	\multicolumn2c{Row likelihood}	&	\multicolumn2c{Row-column likelihood}	\\\hline
10	&	200	&	2	&	2	&	0.5	&	3.317   & (0.362)   & 0.664 & (0.190) & 0.552 & (0.107)	\\
15	&	200	&	2	&	2	&	0.5	&	188.944 & (2.761)   & 0.784 & (0.219) & 0.654 & (0.139)	\\
10	&	400	&	2	&	2	&	0.5	&	5.505   & (0.063)   & 0.784 & (0.154) & 0.805 & (0.122)	\\
10	&	200	&	3	&	2	&	0.5	&	350.976 & (112.915) & 0.501 & (0.153) & 0.539 & (0.113)	\\
10	&	200	&	2	&	3	&	0.5	&	3.441   & (0.672)   & 4.734 & (2.649) & 1.558 & (0.589)	\\
10	&	200	&	2	&	2	&	1.0	&	3.222   & (0.612)   & 1.082 & (0.426) & 0.987 & (0.234)	\\\hline
\end{tabular}}
\caption{Median computing time (and median absolute deviation) in seconds.}\label{tab:simula_time}
\end{table}
We see that  the {\em (row) mass probabilities} $\la_u$ are well estimated by both approximations, with accuracy comparable to full likelihood estimation -- suggesting that, in scenarios with $r<<s$, there is enough information available on each row latent variable for either (and both) of the proposed composite likelihood approximations to accurately capture such probabilities. 

In contrast, the {\em (column) transition probabilities} $\pi_{\bar{v}v}$ are estimated with comparable accuracy by the two approximations, but this accuracy is lower than that afforded by full likelihood estimation -- likely reflecting the fact that, even in the more sophisticated row-column approximation, $c\ell_2(\b\th)$ relies on independent data columns. 

Finally, the {\em means} 
$\psi_{uv}$ are estimated with higher accuracy by the row-column approximation than by the row approximation -- reflecting the fact that the former comes closer to the full likelihood. Similar comments apply to the estimation of $\sigma^2$.

Concerning computing time, our composite likelihood approximations are about 5-fold faster than the full likelihood in the benchmark design. Perhaps most importantly, when we pass to scenarios where $r=15$ (instead of 10) or $k_1=3$ (instead of 2) time increases by two orders of magnitude for the full likelihood. We also note that, while the median running times for the full likelihood appear still relatively modest (approximately $351$ seconds), in some of the simulations with $k_1=3$ they were as high as 9-10 hours -- notwithstanding the fact that size and complexity of the simulated data here are still much smaller than those one can expect in real applications (for an application of the size and complexity of the one in Section~7, running times for the full likelihood could be measured in months). 

This effect of row size and structural complexity is not seen for either of our approximations. Their average computing times remain fairly similar across scenarios, and appear appreciably higher only when $k_2=3$ (instead of 2).

In general, average times for row and row-column approximations are also similar to each other. In fact, in some cases (e.g.~the one with $k_2=3$) the row-column approximation appears to be faster than the row approximation; this is due to the fact that the EM algorithm converges in a smaller number of iterations, even though each iteration is more time consuming by construction.

In summary, our simulations show the row-column composite likelihood approximation to be the right compromise between accuracy and computational viability; it is closer to the accuracy of the full likelihood estimation than the row approximation, especially for estimating means and variance, but much cheaper than the full likelihood for large/complex data arrays -- and not more expensive than the row approximation. 


\section{Model selection}\label{sec:model selection}

A critical point for the model and the composite likelihood approach we propose to be useful in applications, is selecting the number of support points for row ($k_1$) and column ($k_2$) latent variables.
When full likelihood methods are used to estimate simpler models, the literature on model selection suggests information criteria such as the Akaike Information Criterion \citep[AIC;][]{aka:73} and the Bayesian Information Criterion \citep[BIC;][]{sch:78} (see \cite{mclachlan2004finite}, Chapter 6, for a general discussion on selecting the number of components in finite mixture models). These criteria penalize the maximum log-likelihood of the model of interest with a term based on the number of free parameters, seen as a measure of model complexity.

Adaptations of both the AIC and the BIC in which the maximum of the full log-likelihood is replaced with that of a composite log-likelihood are proposed by \cite{var:vid:05} and \citet{gao-song:2011}. In these cases, computing the penalization term is more complicated -- as it requires the Hessian of the composite log-likelihood function and estimation of the variance of its score; see also \cite{bartolucci2015pairwise}.

Given the complexity of the model we introduced, we prefer to rely on a cross validation strategy similar to that in \cite{smyth2000model} and \cite{celeux:durand:2008} -- which avoids the matrices involved in the modified AIC and BIC altogether. In this regard, we note that, since we are not dealing with independent and identically distributed data, estimation of the composite log-likelihood score is rather complicated. On the other hand, cross validation can be implemented straightforwardly, requiring only a small amount of extra code with respect to that already developed for estimation, and a reasonable computing time.

The cross validation strategy we propose, after splitting the data into a training and a validation sample, treats the missing cells in either sample as ``missing completely at random". In more detail, for selecting $k_1$ and $k_2$ we proceed as follows:
\begin{itemize}
\item Split the data into a training sample $\cg S_d$ and a validation sample $\bar{\cg S}_d$ by randomly drawing one half of the cells in the observed two-way array, and repeat this $d=1,\ldots, D$ times (e.g.~$D=100$ is used in our application below). 
\item For each $d=1,\ldots,D$ and each pair $(k_1,k_2)$ of interest, estimate the parameters in $\b\th$ based on $\cg S_d$ by maximizing $c\ell_{k_1 k_2}(\b\th|\cg S_d)$ under the assumption that the cells removed for validation are data missing completely at random. Let $\hat{\b\th}_{k_1k_2}(\cg S_d)$ indicate the resulting estimate.
\item For each pair $(k_1,k_2)$ of interest, compute
\begin{eqnarray*}
c\ell_{cv,k_1k_2} &=& \frac{1}{D}\sum_{d=1}^D c\ell_{k_1k_2}(\hat{\b\th}_{k_1k_2}(\cg S_d)|\bar{\cg S}_d),\\
n_{cv,k_1k_2} &=& \sum_{d=1}^{D} 1\bigg\{c\ell_{k_1k_2}(\hat{\b\th}_{k_1k_2}(\cg S_d)|\bar{\cg S}_d)=
\max_{h_1,h_2}c\ell_{h_1h_2}(\hat{\b\th}_{h_1h_2}(\cg S_d)|\bar{\cg S}_d)\bigg\},
\end{eqnarray*}
where $1\{\cdot\}$ is the indicator function equal to $1$ if its argument is true and to $0$ otherwise. The first quantity
is the average composite log-likelihood computed on the validation samples -- considering for each the parameter estimates based on the corresponding training sample. The second quantity
is 
the number of validation samples (out of $D$) for which the model with $k_1$ and $k_2$ support points reaches the highest value of the composite log-likelihood.
\end{itemize}
As we illustrate in the application section below, these quantities provide guidance in choosing $(k_1,k_2)$; we would like a pair that either maximizes or reaches a value close to the maximum in terms of both $c\ell_{cv,k_1k_2}$ and $n_{cv,k_1k_2}$. Of course other derived quantities, as well as parsimony considerations, can and should be employed also (see below).


\section{A first application to genomic data}\label{sec:illustration}

As a first illustration of how our model and methodology can be used on large, complex data sets, we consider an application to Genomics. The data comes from a study by Kuruppumullage Don et al. (2013) and has been kindly provided by K.D.~Makova and her group at the Pennsylvania State University. The authors used standard HM methodology to segment the human genome based on the rates of four types of mutations estimated from primate comparisons in contiguous 1Mb non-overlapping windows. To try and relate the resulting ``mutational states" to the landscape of DNA composition and molecular activity along the genome, the authors also gathered and pre-processed publicly available data on several dozens genomic features in the same windows system. Here, we address the question of whether it is possible to produce another segmentation, based not on four mutation rates but on this large array of features -- while simultaneously characterizing their interdependencies through clustering. As a feasibility proof, we thus utilize our model and methodology to perform ``clustering-by-segmentation" on a two-way data array comprising $r=28$ features measured in $s=224$ contiguous 1Mb non-overlapping windows covering human chromosome 1. 

\begin{figure}[ht!]\centering
\includegraphics[width=18cm]{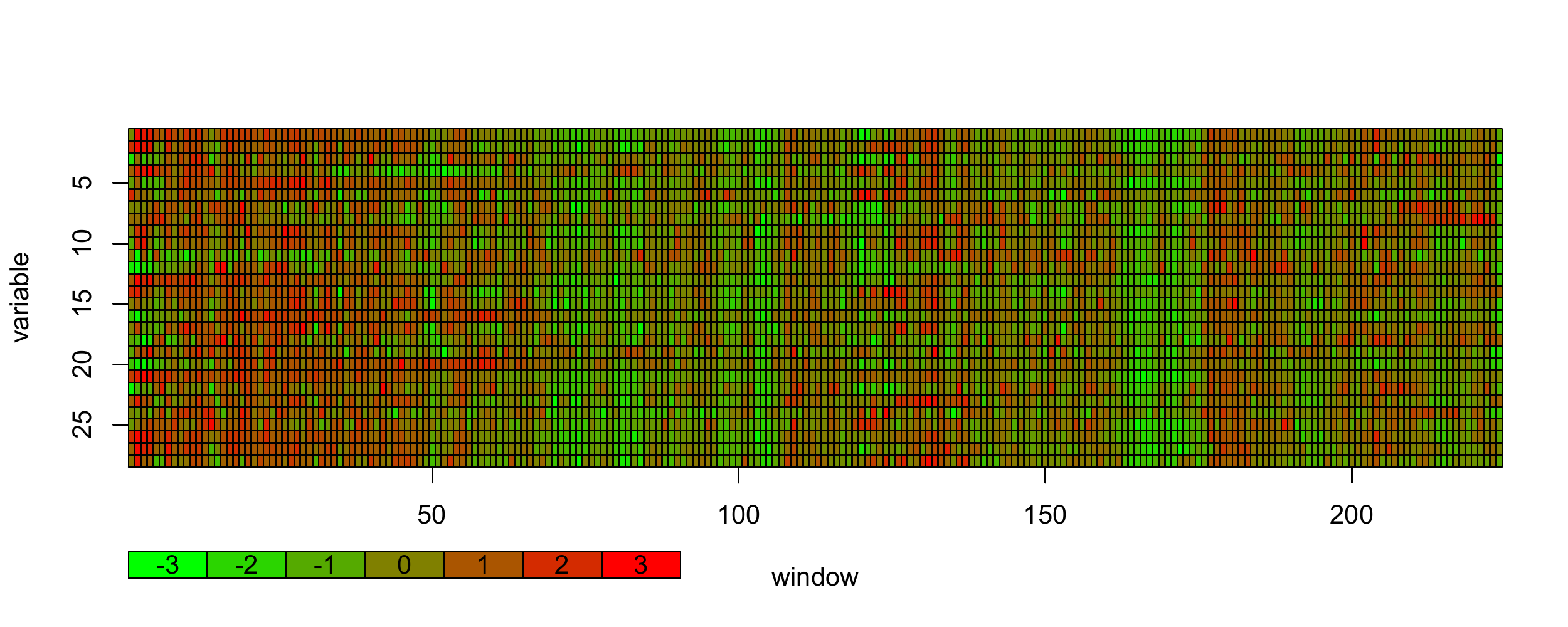}
\caption{\em Data on $r=28$ features measured in $s=224$ contiguous 1Mb non-overlapping windows covering human chromosome 1, after standardization.}\label{fig:data}
\end{figure}

The features, listed in Table~\ref{tab:list_variables}, capture aspects of DNA composition (e.g.~GC content), prevalence of transposable elements (e.g.~number of LINE elements; SINE elements; DNA transposons -- as well as their subfamilies), 
recombination (male and female recombination rates), chromatin structure (e.g.~number of nuclear lamina associated regions; miRNAs; H3K4me1 sites and H3K14 acetylation sites; Polymerase II binding sites; DNase 1 hypersensitive sites), methylation (e.g.~number of non-CpG methylated cytosines; 5-hydroxymethylcitosines; average DNA methylation level), transcription (e.g.~number of CpG islands;  coverage by coding exons) and more. The features were standardized through normal scores prior to use with our approach. A representation of the data after standardization is provided in Figure~\ref{fig:data}. 

\begin{table}[ht!]\centering
{\scriptsize
\begin{tabular}{rllcrrrr}
\hline\hline
\multicolumn1c{$i$}	&	\multicolumn1c{Feature Name} & \multicolumn1c{Description} \\\hline
1	&	GC		&	GC content	\\
2	&	CpG		&	N. CpG islands\\
3	&	nCGm	&	N. non-CpG methyl-cytosines\\
4	&	LINE	&	N. LINE elements\\
5	&	SINE	&	N. SINE elements\\
6	&	NLp		&	N. nuclear lamina associated regions\\
7	&	fRec	&	Female recombination rates	\\
8	&	mRec	&	Male recombination rates	\\
9 	&	H3K4me1	&	N. H3K4me1 sites	\\
10	&	pol2	&	N. RNA polymerase-II binding sites	\\
11	&	telomerase$\_$hex	&	N. telomerase containing hexamers	\\
12	&	dna$\_$trans	&	N. DNA transposons	\\
13	&	X5hMc	&	N. 5-hydroxymethylcytosines	\\
14	&	meth$\_$level	&	Average value of DNA methylation level	\\
15	&	RepT	&	Replication timing in human ES cells\\
16	&	mir	&	N. mammalian interspersed repeat elements (subset of SINEs)	\\
17	&	alu	&	N. Alu elements (subset of SINEs)		\\
18	&	mer	&	N. mammalian dna transposons (subset of dna$\_$trans)	\\
19	&	l1	&	N. L1-elements (subset of LINEs)	\\
20	&	l2	&	N. L2-elements (subset of LINEs)	\\
21	&	l1target	&	N. L1 target sites	\\
22	&	h3k14ac	&	N. Histone H3K14 acetylation sites	\\
23	&	miRNA	&	N. miRNA sites	\\
24	&	triplex	&	N. triplex motifs	\\
25	&	inverted	&	N. inverted repeats	\\
26	&	gquadraplex	&	N. G-Quadruplex structure forming motifs	\\
27	&	dnase1	&	N. dnase-1 hypersensitive sites (from ENCODE. ES cells)	\\
28	&	cExon	&	Coverage by coding exons	\\
\hline
\end{tabular}}
\caption{\em Features in the genomics data set provided by K.D.~Makova and her group at the Pennsylvania State University (see also Kuruppumullage Don et al., 2013).}\label{tab:list_variables}
\end{table}

\subsection{Model selection}
The first critical task is to select the number of support points for the row and column latent variables distributions ($k_1$ and $k_2$); that is, the number of groups in which to cluster the $r=28$ genomic features under consideration, and the number of distinct states in which to segment the $s=224$ windows covering chromosome 1. To perform this selection, we relied on
the cross validation strategy described in Section~6 with $D=100$ iterations; 
Table~\ref{tab:cross1} reports the average row-column composite log-likelihood computed on the validation samples, using the estimates computed on the corresponding training samples; this 
is denoted by $c\ell_{cv,k_1k_2}$.
Table~\ref{tab:cross2} reports the number of times a certain model (i.e.~combination of $k_1$ and $k_2$) beat all other models in terms of composite log-likelihood on the validation samples, denoted by $n_{cv,k_1k_2}$. We also report the number of free parameters for each model in Table~\ref{tab:n_free}.

\begin{table}[ht]
\centering
{\scriptsize
\hspace*{-0.4cm}\begin{tabular}{rrrrrrrrrrrrrrrr}\hline\hline
 & \multicolumn{14}c{$k_2$}\\ \cline{2-15}
$k_1$ &\multicolumn1c{2} &\multicolumn1c{3} & \multicolumn1c{4} & \multicolumn1c{5} & \multicolumn1c{6} & \multicolumn1c{7} & \multicolumn1c{8} & \multicolumn1c{9} & \multicolumn1c{10} & \multicolumn1c{11} & \multicolumn1c{12} & \multicolumn1c{13} & \multicolumn1c{14} & \multicolumn1c{15}\\\hline
1	&	-8357.2	&	-8181.6	&	-8141.7	&	-8131.2	&	-8128.8	&	-8127.8	&	-8126.7	&	-8125.8	&	-8125.3	&	-8125.2	&	-8124.3	&	-8123.8	&	-8124.0	&	-8123.8	\\
2	&	-8342.8	&	-8099.8	&	-8068.5	&	-8046.0	&	-8030.2	&	-8017.5	&	-8008.7	&	-8000.5	&	-7997.2	&	-7992.8	&	-7991.3	&	-7989.5	&	-7988.1	&	-7988.1	\\
3	&	-8328.5	&	-8096.4	&	-8058.2	&	-8019.3	&	-8000.3	&	-7986.5	&	-7979.7	&	-7970.3	&	-7970.5	&	-7969.1	&	\cb-7968.1	&	-7968.9	&	-7969.9	&	
-7969.0	\\
4	&	-8327.7	&	-8095.4	&	-8053.2	&	-8010.4	&	-7993.1	&	-7984.4	&	-7977.0	&	-7971.7	&	-7969.7	&	-7970.8	&	-7972.6	&	-7973.2	&	-7975.1	&	-7977.4	\\
5	&	-8327.5	&	-8094.8	&	-8051.7	&	-8005.5	&	-7993.9	&	-7981.7	&	-7975.7	&	-7974.0	&	-7972.9	&	-7974.5	&	-7976.9	&	-7979.4	&	-7983.1	&	-7988.1	\\
6	&	-8327.4	&	-8094.7	&	-8049.7	&	-8005.3	&	-7989.2	&	-7980.9	&	-7977.3	&	-7974.5	&	-7975.5	&	-7975.3	&	-7979.3	&	-7985.8	&	-7989.3	&	-7995.1	\\
7	&	-8327.0	&	-8093.4	&	-8050.0	&	-8004.2	&	-7990.4	&	-7982.3	&	-7976.3	&	-7973.6	&	-7978.0	&	-7982.0	&	-7987.0	&	-7993.6	&	-7999.8	&	-8006.7	\\
8	&	-8326.5	&	-8093.4	&	-8049.6	&	-8002.6	&	-7990.6	&	-7982.4	&	-7977.9	&	-7981.2	&	-7984.1	&	-7986.4	&	-7994.3	&	-8002.8	&	-8009.2	&	-8020.5	\\
9	&	-8326.7	&	-8093.5	&	-8048.1	&	-8002.7	&	-7990.8	&	-7984.5	&	-7982.8	&	-7982.1	&	-7986.0	&	-7990.4	&	-7998.0	&	-8005.4	&	-8014.5	&	-8030.5	\\
10	&	-8326.0	&	-8093.1	&	-8048.7	&	-8001.1	&	-7989.9	&	-7986.3	&	-7982.1	&	-7985.1	&	-7992.0	&	-7995.6	&	-8005.6	&	-8012.5	&	-8028.7	&	-8040.4	\\
   \hline
\end{tabular}}
\caption{\em Average row-column composite log-likelihood on the validation samples for each $k_1$ and $k_2$ combination, obtained from $D=100$ cross validation replicates. The highest value (highlighted) is achieved for $k_1=3$, $k_2=12$  (for $k_1=k_2=1$, the average composite log-likelihood is $-8887.8$).}\label{tab:cross1}
\end{table}

\begin{table}[ht]
\centering
{\scriptsize
\begin{tabular}{rrrrrrrrrrrrrrrr}
  \hline\hline
  & \multicolumn{14}c{$k_2$}\\\cline{2-15}
$k_1$ & 2 & 3 & 4 & 5 & 6 & 7 & 8 & 9 & 10 & 11 & 12 & 13 & 14 & 15 \\ 
  \hline
1  & 0	& 0	& 0	& 0	& 0	& 0	& 0	& 0	& 0	& 0	& 0	& 0	& 0	& 0\\
2  & 0	& 0	& 0	& 0	& 0	& 0	& 0	& 0	& 0	& 0	& 0	& 1	& 1	& 0\\
3  & 0	& 0	& 0	& 0	& 0	& 0	& 0	& 4	& 2	& 6	& 5	& 5	& 4	& 4\\
4  & 0	& 0	& 0	& 0	& 0	& 1	& 2	& 2	& 1	& 5	& 3	& 3	& 3	& 4\\
5  & 0	& 0	& 0	& 0	& 0	& 0	& 2	& 3	& 1	& 0	& 2	& 6	& 0	& 3\\
6  & 0	& 0	& 0	& 0	& 0	& 1	& 1	& 2	& 3	& 5	& 0	& 1	& 0	& 0\\
7  & 0	& 0	& 0	& 0	& 0	& 1	& 0	& 3	& 1	& 0	& 1	& 0	& 0	& 0\\
8  & 0	& 0	& 0	& 0	& 0	& 0	& 0	& 2	& 0	& 0	& 2	& 0	& 0	& 0\\
9  & 0	& 0	& 0	& 0	& 0	& 1	& 1	& 1	& 1	& 0	& 0	& 0	& 0	& 0\\
10 & 0	& 0	& 0	& 0	& 0	& 0	& 0	& 0	& 0	& 0	& 0	& 0	& 0	& 0\\
   \hline
\end{tabular}}
\caption{\em Number of times (out of 100) in which the model has the highest cross validation composite log-likelihood for each $k_1$ and $k_2$ combination.}\label{tab:cross2}
\end{table}

\begin{table}[ht]
\centering
{\scriptsize
\begin{tabular}{rrrrrrrrrrrrrrrr}
  \hline\hline
  & \multicolumn{14}c{$k_2$}\\\cline{2-15}
$k_1$ & 2 & 3 & 4 & 5 & 6 & 7 & 8 & 9 & 10 & 11 & 12 & 13 & 14 & 15 \\ 
  \hline
1	&	5	&	10	&	17	&	26	&	37	&	50	&	65	&	82	&	101	&	122	&	145	&	170	&	197	&	226	\\
2	&	8	&	14	&	22	&	32	&	44	&	58	&	74	&	92	&	112	&	134	&	158	&	184	&	212	&	242	\\
3	&	11	&	18	&	27	&	38	&	51	&	66	&	83	&	102	&	123	&	146	&	171	&	198	&	227	&	258	\\
4	&	14	&	22	&	32	&	44	&	58	&	74	&	92	&	112	&	134	&	158	&	184	&	212	&	242	&	274	\\
5	&	17	&	26	&	37	&	50	&	65	&	82	&	101	&	122	&	145	&	170	&	197	&	226	&	257	&	290	\\
6	&	20	&	30	&	42	&	56	&	72	&	90	&	110	&	132	&	156	&	182	&	210	&	240	&	272	&	306	\\
7	&	23	&	34	&	47	&	62	&	79	&	98	&	119	&	142	&	167	&	194	&	223	&	254	&	287	&	322	\\
8	&	26	&	38	&	52	&	68	&	86	&	106	&	128	&	152	&	178	&	206	&	236	&	268	&	302	&	338	\\
9	&	29	&	42	&	57	&	74	&	93	&	114	&	137	&	162	&	189	&	218	&	249	&	282	&	317	&	354	\\
10	&	32	&	46	&	62	&	80	&	100	&	122	&	146	&	172	&	200	&	230	&	262	&	296	&	332	&	370	\\
   \hline
\end{tabular}}
\caption{\em Number of free parameters for each $k_1$ and $k_2$ combination.}\label{tab:n_free}
\end{table}

According to these results, the model achieving highest average composite log-likelihood, is the one with $k_1=3$ and $k_2=12$. This model does also well by beating all other models $5$ times (out of $D=100$) 
-- the maximum 
here is 6, which is obtained for $k_1=3$, $k_2=11$ and $k_1=5$ and $k_2=13$. However, from both Table~\ref{tab:cross1} and Table~\ref{tab:cross2} we can see that several alternative $(k_1, k_2)$ pairs provide very similar performance. In addition, from Table~\ref{tab:n_free} we can see that the model with $k_1=3$ and $k_2=12$ has a very large number of free parameters compared to other models with similar performance. To provide an alternative quantification, for each model we compute an index of relative performance. In more detail, for every given combination of $k_1$ and $k_2$ we consider the average composite log-likelihood across cross-validation iterations, subtract the minimum of such quantity over all combinations considered, and divide by the difference between its maximum and minimum:
\[
q_{k_1k_2} = \frac{c\ell_{cv,k_1k_2}-\min_{h_1h_2}c\ell_{cv,h_1h_2}}
{\max_{h_1h_2}c\ell_{cv,h_1h_2}-\min_{h_1h_2}c\ell_{cv,h_1h_2}};
\]
the higher this index, the better the model identified by $k_1$ and $k_2$. Table~\ref{tab:cross3} reports the index values. 

\begin{table}[ht]
\centering
{\scriptsize
\begin{tabular}{rrrrrrrrrrrrrrr}
\hline\hline
 & \multicolumn{14}c{$k_2$}\\ \cline{2-15}
$k_1$ & \multicolumn1c{2} & \multicolumn1c{3} & \multicolumn1c{4} & \multicolumn1c{5} & \multicolumn1c{6} & \multicolumn1c{7} & \multicolumn1c{8} & \multicolumn1c{9} & \multicolumn1c{10} & \multicolumn1c{11} & \multicolumn1c{12} & \multicolumn1c{13} & \multicolumn1c{14} & \multicolumn1c{15} \\\hline
1	&	0.577	&	0.768	&	0.811	&	0.823	&	0.825	&	0.826	&	0.828	&	0.828	&	0.829	&	0.829	&	0.830	&	0.831	&	0.830	&	0.831	\\
2	&	0.593	&	0.857	&	0.891	&	\cb0.915	&	\cb0.932	&	\cb0.946	&	\cb0.956	&	\cb0.965	&	\cb0.968	&	\cb0.973	&	\cb0.975	&	\cb0.977	&	\cb0.978	&	\cb0.978	\\
3	&	0.608	&	0.860	&	\cb0.902	&	\cb0.944	&	\cb0.965	&	\cb0.980	&	\cb0.987	&	\cb0.998	&	\cb0.997	&	\cb0.999	&	\cb1.000	&	\cb0.999	&	\cb0.998	&	\cb0.999	\\
4	&	0.609	&	0.862	&	\cb0.907	&	\cb0.954	&	\cb0.973	&	\cb0.982	&	\cb0.990	&	\cb0.996	&	\cb0.998	&	\cb0.997	&	\cb0.995	&	\cb0.994	&	\cb0.992	&	\cb0.990	\\
5	&	0.609	&	0.862	&	\cb0.909	&	\cb0.959	&	\cb0.972	&	\cb0.985	&	\cb0.992	&	\cb0.994	&	\cb0.995	&	\cb0.993	&	\cb0.990	&	\cb0.988	&	\cb0.984	&	\cb0.978	\\
6	&	0.609	&	0.862	&	\cb0.911	&	\cb0.960	&	\cb0.977	&	\cb0.986	&	\cb0.990	&	\cb0.993	&	\cb0.992	&	\cb0.992	&	\cb0.988	&	\cb0.981	&	\cb0.977	&	\cb0.971	\\
7	&	0.610	&	0.864	&	\cb0.911	&	\cb0.961	&	\cb0.976	&	\cb0.984	&	\cb0.991	&	\cb0.994	&	\cb0.989	&	\cb0.985	&	\cb0.979	&	\cb0.972	&	\cb0.966	&	\cb0.958	\\
8	&	0.610	&	0.864	&	\cb0.911	&	\cb0.962	&	\cb0.975	&	\cb0.984	&	\cb0.989	&	\cb0.986	&	\cb0.983	&	\cb0.980	&	\cb0.971	&	\cb0.962	&	\cb0.955	&	\cb0.943	\\
9	&	0.610	&	0.864	&	\cb0.913	&	\cb0.962	&	\cb0.975	&	\cb0.982	&	\cb0.984	&	\cb0.985	&	\cb0.980	&	\cb0.976	&	\cb0.967	&	\cb0.959	&	\cb0.949	&	\cb0.932	\\
10	&	0.611	&	0.864	&	\cb0.912	&	\cb0.964	&	\cb0.976	&	\cb0.980	&	\cb0.985	&	\cb0.981	&	\cb0.974	&	\cb0.970	&	\cb0.959	&	\cb0.952	&	\cb0.934	&	\cb0.921	\\   \hline
\end{tabular}}
\caption{\em Relative 
performance index for each $k_1$ and $k_2$ combination. Values $\geq 0.9$ (highlighted) are already achieved using fairly few row and column support points.}\label{tab:cross3}
\end{table}

The relative performance index points towards the model with $k_1=3$ and $k_2=4$.
This model achieves $q_{3,4}=0.902$ (i.e.~a loss of predictive power of only 10\% relative to the model with $k_1=3$ and $k_2=12$) while requiring only $27$ free parameters (compared to $171$ for $k_1=3$ and $k_2=12$). In fact, the model with $k_1=3$ and $k_2=4$ is the smallest with a relative performance above $0.9$. Based on cross validation performance and parsimony, we therefore take this as our selected model.

\subsection{Estimation results}
Next, we discuss parameter estimates for our selected model; recall that we are forming three clusters of genomic features ($k_1=3$), and segmenting chromosome 1 according to four distinct states ($k_2=4)$.

Table~\ref{tab:est1} reports estimates of the mass probabilities of the row latent variable distribution ($\hat{\la}_u$) and estimates of the means ($\hat{\psi}_{uv}$). As a convention, 
modalities of the row latent variable ($u=1,2,3$) are ordered by decreasing $\hat{\la}_u$ 
and modalities of the column latent variable ($v=1, 2, 3, 4$) are ordered by increasing $\hat{\psi}_{1v}$.

\begin{table}[ht]
\centering
{\scriptsize
\begin{tabular}{rrrrrrrrrrrrrr}
  \hline\hline
  \vspace*{-2mm}\\
 $u$ & \multicolumn1c{$\hat{\la}_u$} & \multicolumn1c{$\hat{\psi}_{u1}$} & \multicolumn1c{$\hat{\psi}_{u2}$} & \multicolumn1c{$\hat{\psi}_{u3}$} & \multicolumn1c{$\hat{\psi}_{u4}$} 
 \\ \hline
1 & 0.788 & -1.383 & -0.492 & 0.134 & 0.995 \\ 
  2 & 0.132 & -0.574 & 0.832 & 0.101 & -1.022 \\ 
  3 & 0.080 & 1.015 & -0.954 & 0.091 & 0.421 \\
\hline 
\end{tabular}}
\caption{\em Estimates of the mass probabilities of the row latent variable distribution and estimates of the means for the selected model ($k_1=3$ and $k_2=4$). 
}\label{tab:est1}
\end{table}

Table~\ref{tab:est2} reports estimates of the transition probabilities ($\hat{\pi}_{\tilde{v}v}$) and estimates of the stationary distribution ($\hat{\rho}_v$) for the Markov process governing the column latent variable.  

\begin{table}[ht]
\centering
{\scriptsize
\begin{tabular}{rrrrrrrrrrrrrr}
  \hline\hline
$v$ & \multicolumn1c{$\hat{\rho}_v$} & \multicolumn1c{$\hat{\pi}_{v1}$} & \multicolumn1c{$\hat{\pi}_{v2}$} & \multicolumn1c{$\hat{\pi}_{v3}$} & \multicolumn1c{$\hat{\pi}_{v4}$}
\\ \hline
 1 & 0.121 & 0.802 & 0.196 & 0.002 & 0.000 \\ 
  2 & 0.285 & 0.085 & 0.713 & 0.163 & 0.039 \\ 
  3 & 0.339 & 0.000 & 0.171 & 0.797 & 0.032 \\ 
  4 & 0.255 & 0.000 & 0.000 & 0.086 & 0.914 \\\hline
\end{tabular}}
\caption{\em Estimates of the transition probabilities and estimates of the stationary distribution for the Markov process governing the column latent variable for the selected model ($k_1=3$ and $k_2=4$). 
}\label{tab:est2}
\end{table}

\begin{figure}[ht!]\centering
\includegraphics[width=18cm]{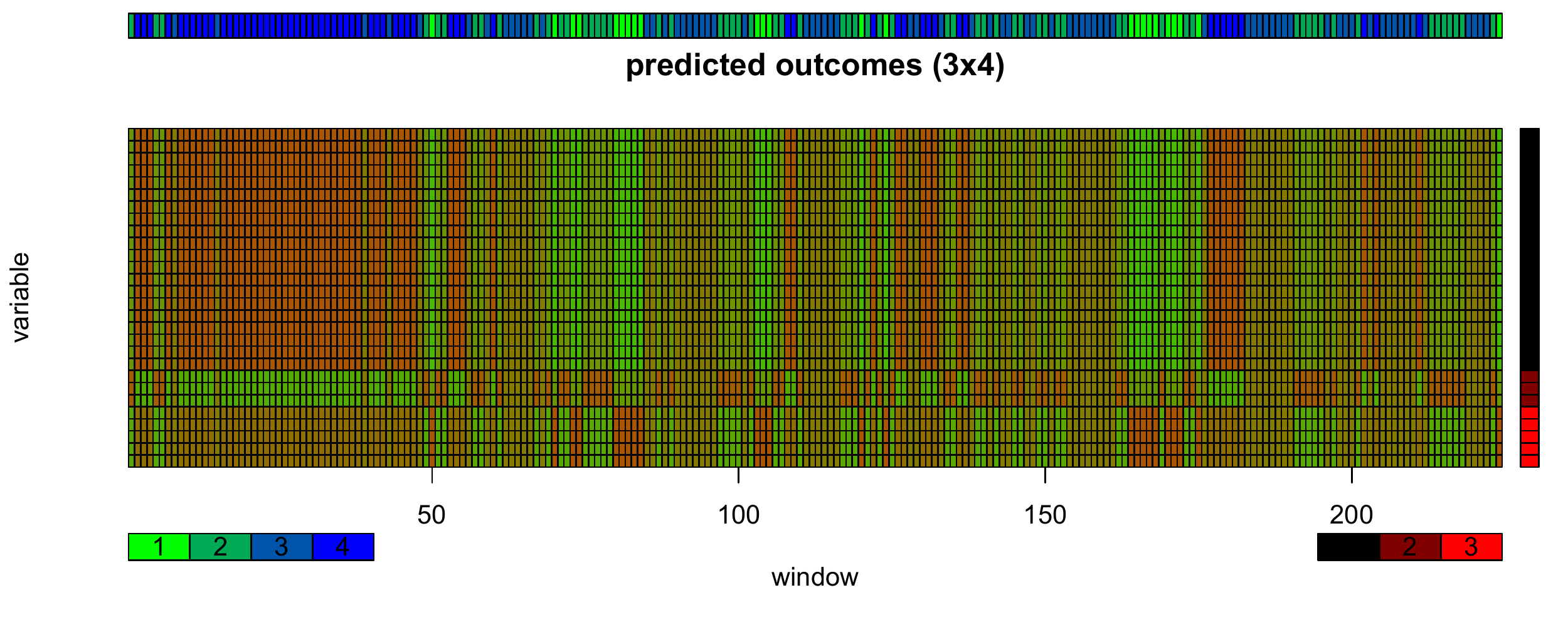}
\caption{\em Color-coded map of predicted genomic feature clusters (right), segmentation states for the windows along chromosome 1 (top) and means of each feature in each window (middle) for the selected model ($k_1=3$ and $k_2=4$). Rows are rearranged according to the assigned clusters -- Cluster 2 comprises number of telomerase containing examers, DNA transposons and histone H3K14 acetylation sites; Cluster 3 comprises number of non-CpG methyl-cytosines, nuclear lamina regions, polymerase II binding sites, ALU elements and MER elements (Cluster 1 groups all the remaining 20 features).
}\label{fig:predict3}
\end{figure}

Figure~\ref{fig:predict3} shows a color-coded map of the predictions associated with the selected model. For each cell $(i,j)$, $i=1,\ldots,r$ ($r=28$), $j=1,\ldots,s$ ($s=224$) of the two-way array, we (i) predict the feature cluster (i.e.~the row latent state) $\hat{u}_i$ and the segmentation state (i.e.~column latent state) $\hat{v}_j$ on the basis of the {\em maximum a posteriori probability} (MAP), and (ii) set the cell's predicted value to the estimated mean $\hat{\psi}_{\hat{u}_i\hat{v}_j}$. The horizontal dimension represents the $s=224$ contiguous windows along chromosome 1, with the horizontal bar on top reporting $\hat{v}_j$'s color-coded on a green-to-blue range. The vertical dimension represents the $r=28$ genomic features, with the vertical bar on the right reporting $\hat{u}_i$'s color-coded on a black-to-red range. Rows are rearranged grouping features according to the three clusters. The inner part of the figure reports the $\hat{\psi}_{\hat{u}_i\hat{v}_j}$'s color-coded on a green-to-red range as was done for the data in Figure~\ref{fig:data}. One can therefore interpret patterns in the way low (green) and high (red) predicted values characterize different genomic feature clusters (as marked on the vertical bar to the right) and segments on the chromosome (as marked on the horizontal bar on top). 

Concerning the three clusters of genomic features, we note that Cluster 1 is very large, comprising 20 features (the estimated mass probability is approximately 80\%), while Clusters 2 and 3 are much smaller, with 3 and 5 features respectively (the estimated mass probabilities are each approximately 10\%). In more detail, Cluster 2 includes number of telomerase containing examers (a proxy for repair), DNA transposons (a proxy for transposition activity) and histone H3K14 acetylation sites (a proxy for chromatin structure). Cluster 3 includes number of non-CpG methyl-cytosines (a proxy for methylation), nuclear lamina regions and polymerase II binding sites (proxies for chromatin structure), ALU elements and MER elements (proxies for transposition activity).

Concerning the four segmentation states, we note that they cover approximately 10\%, 30\%, 35\% and 25\% of chromosome 1, respectively (from the estimated stationary distribution). From the estimated means, we note that State 1, i.e.~the least prevalent, is characterized by strongly depressed Cluster 1 features, depressed Cluster 2 features and strongly elevated Cluster 3 features. State 3, i.e.~the most prevalent, has mildly elevated levels for features in all clusters. State 2 and State 4, whose prevalences are more similar, have ``mirroring" profiles -- the former is characterized by depressed Cluster 1 features, strongly elevated Cluster 2 features and strongly depressed Cluster 3 features, the latter by strongly elevated Cluster 1 features, strongly depressed Cluster 2 features and elevated Cluster 3 features.

Interestingly, from Figure~\ref{fig:predict3} we can see that while all four states are represented and alternate along most of the chromosome, its ``beginning" (approximately the first $50$ windows towards the left of the figure) shows a marked prevalence of State 4. Also interestingly, Cluster 2 shows strongly elevated levels in State 2 (covering approximately 30\% of the chromosome), where all other features are depressed or strongly depressed, and strongly depressed levels in State 4 (covering approximately 25\% of the chromosome with a prevalence in the first $~50$ windows), where all other features are elevated or strongly elevated. On its end, Cluster 3 shows strongly elevated levels in State 1 (covering approximately 10\% of the chromosome), where all other features are depressed or strongly depressed. 




\section{Conclusions}\label{sec:conclusions}

In this article, we considered a discrete latent variable model for two-way data arrays, which allows one to simultaneously produce clusters along one of the data dimensions and contiguous groups, or segments, along the other. We proposed two composite likelihood approximations and their EM-based optimization for estimation, as well as a specialized cross validation strategy to select the number of support points for row and column latent variables.

Through simulations, we showed that our composite likelihood methodology has reasonable performance in comparison with full likelihood methodology (when the latter is viable) while being much less computationally demanding. Our simulations also demonstrated a clear advantage of the row-column composite likelihood with respect to the row (only) composite likelihood in terms of estimation efficiency -- and sometimes also in terms of computing time.

Importantly, our methodology remains computationally viable even when, due the dimension or the structural complexity of the data, the full likelihood cannot be used; this is likely to happen in many practical applications -- especially in ones involving genomic data, such as the one we presented in Section~\ref{sec:illustration}.

Another important consequence of the low computational burden of our approach is that repeated estimation, such as that required in cross validation, may be run in a reasonable computing time. This allowed us to implement model selection using a straightforward cross validation strategy.

Our first application to genomic data, albeit preliminary, demonstrated the feasibility of using composite likelihood methodology to simultaneously segment long stretches of a genome and cluster large arrays of genomic features. For instance, we were able to identify about $50$Mb at the beginning of human chromosome 1 where most of the $28$ genomic features we considered tend to be elevated or strongly elevated, but three (number of telomerase containing examers, a proxy for repair; DNA transposons, a proxy for transposition activity; and histone H3K14 acetylation sites, a proxy for chromatin structure) tend to be strongly depressed. A similar analysis could be extended to all chromosomes and a yet broader set of features, unveiling important biological clues.

Our model and methodology could also be used on many other types of complex genomic data, and applied to many other fields. For instance, they could be used for analyzing parallel time-series of economic indicators recorded on several countries (see Introduction), or data from Item Response Theory (rows corresponding to examinees, columns corresponding to test items administered sequentially).

Regarding further methodological developments, we plan to explore more sophisticated forms of composite likelihood approximation, which may lead to additional improvements in estimation efficiency -- e.g.~in estimating transition probabilities for the Markov process governing the column latent variable, for which estimation quality with our row-column and row composite likelihood appeared poorer than for other parameters in simulations. 


\section*{Acknowledgements}
We are grateful to K.D.~Makova (the Pennsylvania State University) and her group, who provided data and help for the Genomics application presented in this article. Our ongoing collaboration creates a rich and motivating context for methodological research. F.~Bartolucci acknowledges financial support from award RBFR12SHVV of the Italian Government (FIRB ``Mixture and latent variable models for causal inference and analysis of socio-economic data", 2012). F.~Chiaromonte and P.~Kuruppumullage Don were partially supported by award DBI-0965596 of the U.S.~National Science Foundation (ABI ``Computational tools and statistical analysis of co-varying rates of different mutation types", 2010) and by funds from the Huck Institutes of the Life Sciences of the Pennsylvania State University. B.G.~Lindsay and P.~Kuruppumullage Don were partially supported by funds from the Willaman Professorship and the Eberly Family Chair in Statistics at the Pennsylvania State University.

During the final stages of preparation of this article, B.G.~Lindsay passed away due to an illness. We lost a dear friend, a generous mentor and a brilliant colleague whose insight, rigor and love for our discipline we would like to honor -- and will forever treasure. Bruce's contributions to Statistics and to the lives of so many around him have been invaluable; we will deeply miss him.   
%


\bibliography{biblio_two-way_lvm}
\bibliographystyle{apalike}

\end{document}